\documentclass[]{article}

\usepackage{amssymb}
\usepackage{amsmath}

\title{A general proof system for logics of imperfect information}
\author{Pietro Galliani\\ILLC\\University of Amsterdam\\The Netherlands\\(pgallian{@}gmail.com)}

\newtheorem{theorem}{Theorem}[section]
\newtheorem{definition}[theorem]{Definition}
\newtheorem{lemma}[theorem]{Lemma}
\newtheorem{proposition}[theorem]{Proposition}
\newtheorem{corollary}[theorem]{Corollary}

\newenvironment{proof}{\emph{Proof:}$\\$}{$\\\Box\\$}

\newcommand{\tuple}[1]{\vec{#1}}
\newcommand {\indep}[3] {#2 ~\bot_{#1}~ #3}

\def\dom{\mbox{Dom}}
\def\free{\mbox{Free}}
\def\rel{\mbox{Rel}}
\def\part{\mathcal{P}}

\def\nnf{\mbox{NNF}}

\def\gts{\mbox{GTS}}
\def\lts{\mbox{LTS}}
\def\ens{\mbox{ENS}}

\begin{document}
\maketitle
\begin{abstract}
We develop a semantics for logics of imperfect information with respect to general models.

Then we build a proof system and prove its soundness and completeness with respect to this semantics.
\end{abstract}
\section{Introduction}
Logics of imperfect information are extensions of first-order logic (or, sometimes, of other logics: see for example Tulenheimo \cite{tulenheimo04} and V\"a\"an\"anen \cite{vaananen08b}) which allow to reason about patterns of dependence and independence between variables.\\

Historically, the earliest such logic was \emph{branching quantifier logic} (Henkin \cite{henkin61}), which adds to the language of first order logic branching quantifiers such as 
\[
	\left(\begin{array}{c c}
		\forall x & \exists y\\
		\forall z & \exists w
	\end{array}
	\right)
	\phi(x,y,z,w)
\]
whose interpretation, informally speaking, states that the choice of $y$ is not dependent on the choice of $z$ and the choice of $w$ is not dependent on the choice of $x$. A significant breakthrough in the study of this class of logics occurred with the development of \emph{independence-friendly logic} (Hintikka and Sandu \cite{hintikkasandu89}), through which
\begin{enumerate}
\item The syntax of branching quantifier logic was significantly simplified, doing away with complex structures of quantifiers such as the above one and introducing instead \emph{slashed quantifiers} $(\exists x / W) \phi$, whose informal interpretation is ``there exists a $x$, not dependent on any variables in $W$, such that $\phi$''; 
\item The \emph{game-theoretic semantics} of logics of imperfect information was defined formally, and its properties were examined in detail. 
\end{enumerate}
These developments made it possible to define, in \cite{hodges97}, a \emph{compositional semantics} for independence-friendly logic which is equivalent to its game-theoretic semantics. This semantics, called \emph{team semantics} or \emph{trump semantics}, differs from Tarski's semantics for first order logic in that satisfaction conditions of formulas are predicated not of single assignments, but of \emph{sets} of assignments\footnote{Later, Cameron and Hodges \cite{cameron01} proved, through combinatorial methods, that no compositional semantics for such a logic exists in which the satisfaction conditions are predicated of single assignments.} (which we will henceforth call \emph{Teams}, after the terminology of V\"a\"an\"anen \cite{vaananen07}).\\

This alternate semantics provided one of the main impulses towards the development of \emph{dependence logic} \cite{vaananen07}, which separates the notion of dependence and independence from the notion of quantification by doing away with slashed quantifier and introducing instead \emph{dependence atoms} of the form $=\!\!(t_1 \ldots t_n)$, where $t_1 \ldots t_n$ are terms, which are satisfied by a team $X$ if and only if the value of $t_n$ is a function of the values of $t_1 \ldots t_{n-1}$ in it. This - only at first sight minor - innovation led to a number of significant advances in the study of the properties of logics of imperfect information, and, in particular, of their model theory; apart from the aforementioned \cite{vaananen07}, we can refer here for example to the results of (Juha) Kontinen and V\"a\"an\"anen \cite{kontinenv09} and (Jarmo) Kontinen \cite{kontinen_ja10}.\\

Furthermore, a recent direction of research in the field of logics of imperfect information consists in the study of the model-theoretical properties of variants of dependence logic obtained by substituting the dependence atoms with other kinds of non first-order atomic formulas. The earliest work along these lines was Gr\"adel and V\"a\"an\"anen \cite{gradel10}, whose \emph{independence logic} is expressively stronger than dependence logic and will be the main logical formalism taken in exam in the rest of this work; furthermore, we have \emph{multivalued dependence logic} from Engstr\"om \cite{engstrom10} and \emph{inclusion logic} and \emph{exclusion logic} from Galliani \cite{galliani11}.\footnote{That paper also characterized precisely the expressive power of independence logic with respect to open formulas, thus answering an open problem of \cite{gradel10}, and proved that inclusion and exclusion logic are strictly weaker than independence logic.}\\

One property common to all these papers is that they are essentially concerned only with the \emph{semantics} of logics of imperfect information and its model-theoretic properties. The corresponding proof theories, instead, are still relatively undeveloped. The recent \cite{kontinenv11} presents a sound and complete deduction system for extracting the first-order consequences of a Dependence Logic theory; however, due to the equivalence between Dependence Logic and existential second-order logic there exists no hope of extending this system to one for deducing the Dependence Logic consequences of a Dependence Logic theory according to the standard semantics. The present paper, drawing inspiration from Henkin's treatment of second order logic \cite{henkin50} and from the analysis of branching quantifiers of \cite{lopezescobar91}, may be seen as a different approach to the study of the proof theories of logics of imperfect information: instead of restricting our language, we will weaken the semantics and consider a more general class of models and then we will develop a proof system capable of extracting all valid formula for this new semantics. 
\section{Independence Logic}
In this section, we will briefly recall the syntax and the semantics of Independence Logic, plus a few of its basic properties. It can be safely skipped by anyone who is already familiar with the results of \cite{gradel10}.\\

As is often done in the field of logics of imperfect information, we will assume that our expressions are always in Negation Normal Form.
\begin{definition}[Syntax]
	Let $\Sigma$ be a first order signature. Then the set $\nnf_\Sigma$ of the \emph{negation normal form formulas} of our logic is the smallest set such that 
	\begin{description}
		\item[NNF-lit] If $\phi$ is a first order literal over the signature $\Sigma$ then $\phi \in \nnf_\Sigma$;
		\item[NNF-ind] If $\tuple t_1$, $\tuple t_2$ and $\tuple t_3$ are tuples of terms with signature $\Sigma$ then $\indep{\tuple t_1}{\tuple t_2}{\tuple t_3}$ is in $\nnf_{\Sigma}$;
		\item[NNF-$\vee$] If $\phi$ and $\psi$ are in $\nnf_\Sigma$ then $\phi \vee \psi$ is also in $\nnf_\Sigma$;
		\item[NNF-$\wedge$] If $\phi$ and $\psi$ are in $\nnf_\Sigma$ then $\phi \wedge \psi$ is also in $\nnf_\Sigma$;
		\item[NNF-$\exists$] If $\phi$ is in $\nnf_\Sigma$ and $x$ is a variable then $\exists x \phi$ is in $\nnf_\Sigma$; 
		\item[NNF-$\forall$] If $\phi$ is in $\nnf_\Sigma$ and $x$ is a variable then $\forall x \phi$ is in $\nnf_\Sigma$.
	\end{description}
\end{definition}
The set $\free(\phi)$ of the \emph{free variables} of a formula $\phi$ is defined similarly to the case of First Order Logic:
\begin{definition}[Free Variables]
	Let $\Sigma$ be a first order signature, and let $\phi \in \nnf_\Sigma$. Then the set $\free(\phi)$ of the \emph{free variables} of $\phi$ is defined by structural induction on $\phi$ as follows:
	\begin{description}
		\item[Free-lit] If $\phi$ is a first order literal then $\free(\phi)$ is the set of all variables occurring in $\phi$; 
		\item[Free-ind] If $\phi$ is $\indep{\tuple t_1}{\tuple t_2}{\tuple t_3}$ then $\free(\phi)$ is the set of all variables occurring in $\tuple t_1$, $\tuple t_2$ or $\tuple t_3$;
		\item[Free-$\vee$] If $\phi$ is $\psi \vee \theta$ for some formulas $\psi, \theta \in \nnf_{\Sigma}$ then $\free(\phi)$ is $\free(\psi) \cup \free(\theta)$; 
		\item[Free-$\wedge$] If $\phi$ is $\psi \wedge \theta$ for some formulas $\psi, \theta \in \nnf_{\Sigma}$ then $\free(\phi)$ is $\free(\psi) \cup \free(\theta)$; 
		\item[Free-$\exists$] If $\phi$ is $\exists x \psi$ for some variable $x$ and some $\psi \in \nnf_\Sigma$ then $\free(\phi) = \free(\psi) \backslash \{x\}$;
		\item[Free-$\forall$] If $\phi$ is $\forall x \psi$ for some variable $x$ and some $\psi \in \nnf_\Sigma$ then $\free(\phi) = \free(\psi) \backslash \{x\}$.
	\end{description}
\end{definition}

The following definition is standard: 
\begin{definition}[Team]
	Let $V$ be a finite set of variables and let $M$ be a first order model. A \emph{team} over $M$ with domain $V$ is a set of first order assignments over $M$ with domain $V$. 
\end{definition}

The next definition will be useful to give the semantics for the ``lax'' (in the sense of \cite{galliani11}) version of the existential quantifier that we will use: 
\begin{definition}[$x$-variation]
	Let $M$ be a first order model, let $X$ be a team over $M$, and let $x$ be a variable symbol (not necessarily in $\dom(X)$). Then a team $X'$ of $M$ with domain $\dom(X') = \dom(X) \cup \{x\}$ is said to be a $x$\emph{-variation} of $X$, and we write $X[x]X'$, if and only if the restrictions of $X$ and $X'$ to $\dom(X) \backslash \{x\}$ are the same.
\end{definition}

At this point, we have all that we need in order to define the team semantics for independence logic.
\begin{definition}[Team Semantics for Independence Logic]
	Let $\Sigma$ be a first order signature, let $M$ be a first order model of signature $\Sigma$, let $\phi \in \nnf_\Sigma$ and let $X$ be a team with domain containing $\free(\phi)$. Then we say that $X$ \emph{satisfies} $\phi$ in $M$, and we write $M \models_X \phi$, if and only if 
	\begin{description}
		\item[TS-lit] $\phi$ is a first order literal and, for all $s \in X$, $M \models_s \phi$ in the usual first order sense; 
		\item[TS-ind] $\phi$ is $\indep{\tuple t_1}{\tuple t_2}{\tuple t_3}$ for some tuples of terms $\tuple t_1$, $\tuple t_2$ and $\tuple t_3$, and for all $s, s' \in X$ with $\tuple t_1\langle s\rangle = \tuple t_1\langle s'\rangle$ there exists a $s'' \in X$ with $\tuple t_1 \tuple t_2\langle s''\rangle = \tuple t_1 \tuple t_2\langle s\rangle$ and $\tuple t_1 \tuple t_3\langle s''\rangle = \tuple t_1 \tuple t_3\langle s'\rangle$;
		\item[TS-$\vee$] $\phi$ is $\psi_1 \vee \psi_2$ for two formulas $\psi_1, \psi_2 \in \nnf_{\Sigma}$ and $X = Y \cup Z$ for some two teams $Y$ and $Z$ such that $M \models_Y \psi_1$ and $M \models_Z \psi_2$;
		\item[TS-$\wedge$] $\phi$ is $\psi_1 \wedge \psi_2$ for two formulas $\psi_1, \psi_2 \in \nnf_{\Sigma}$, $M \models_X \psi_1$ and $M \models_X \psi_2$;
		\item[TS-$\exists$] $\phi$ is $\exists x \psi$ for some variable $x$ and some $\psi \in \nnf_\Sigma$ and there exists a team $X'$ such that $X[x]X'$ (that is, $X'$ is a $x$-variation of $X$) and such that $M \models_{X'} \psi$; 
		\item[TS-$\forall$] $\phi$ is $\forall x \psi$ for some suitable $x$ and $M \models_{X[M/x]} \psi$, where 
			\[
			X[M/x] = \{s[m/x] : s \in X, m \in \dom(M)\}.
			\]
	\end{description}
\end{definition}
As \cite{gradel10} shows, the dependence atom $=\!\!(t_1 \ldots t_n)$ is equivalent to the independence atom $\indep{t_1 \ldots t_{n-1}}{t_n}{t_n}$. Therefore, Dependence Logic is contained in Independence Logic. The following result is also in \cite{gradel10}:
\begin{theorem}[\cite{gradel10}]
	Let $\Sigma$ be a first order signature, let $V = \{\tuple v\}$ be a finite set of variables and let $\phi(\tuple v) \in \nnf_\Sigma$ be an Independence Logic formula with signature $\Sigma$ and free variables in $V$. Then there exists an existential second order logic formula $\Phi(R)$ such that, for all models $M$ with signature $\Sigma$ and all teams $X$ over $M$ with domain $V$, 
	\[
		M \models_X \phi \Leftrightarrow M \models \Phi(\rel(X))
	\]
	where $\rel(X) = \{s(\tuple v) : s \in X\}$.
\end{theorem}
In \cite{galliani11}, the converse of this result is proved:
\begin{theorem}[\cite{galliani11}]
	Let $\Sigma$ be a first order signature, let $V = \{\tuple v\}$ be a finite set of variables and let $\Phi(R)$ be an existential second order formula with signature $\Sigma$ and with $R$ as its only free variable, where $R$ is a relational variable of arity $|\tuple v|$. Then there exists an independence logic formula $\phi(\tuple v)$, over the signature $\Sigma$ and with free variables in $\tuple v$, such that 
	\[
		M \models_X \phi \Leftrightarrow M \models \Phi(\rel(X))
	\]
	for all models $M$ with signature $\Sigma$ and all nonempty teams $X$ over $M$ with domain $V$.
\end{theorem}
\section{General models for independence logic}
In this section, we will develop a generalization of team semantics, along the lines of Henkin's treatment of second order logic. As we will see, the fact that Independence Logic corresponds to Existential Second Order Logic (and not to full Second Order Logic) means that we will be able to restrict ourselves to consider only a very specific kind of general model.\\

\begin{definition}[General Model]
	Let $\Sigma$ be a first order signature. A \emph{general model} with signature $\Sigma$ is a pair $(M, \mathcal G)$, where $M$ is a first order model with signature $\Sigma$ and $\mathcal G$ is a set of teams over finite domains, respecting the condition
	\begin{itemize}
		\item If $\phi(x_1 \ldots x_n, \tuple m, \tuple R)$ is a first order formula, where $\tuple m$ is a tuple of constant parameters in $\dom(M)$ and where $\tuple R$ is a tuple of ``relation parameters'' corresponding to teams in $\mathcal G$, in the sense that each $R_i$ is of the form 
			\[
				R_i = \rel(X) = \{s(\tuple z) : s \in X_i\}
			\]
			for some $X_i \in \mathcal G$, then for
			\[
				\|\phi(x_1 \ldots x_n, \tuple m, \tuple R)\|_M = \{s : \dom(s) = \{x_1 \ldots x_n\}, M \models_s \phi(x_1 \ldots x_n, \tuple m, \tuple R)\}
			\]
			it holds that $\|\phi(x_1 \ldots x_n, \tuple m, \tuple R)\|_M \in \mathcal G$. 
	\end{itemize}
\end{definition}

\begin{lemma}
	Let $\Sigma$ be a first order signature and let $(M, \mathcal G)$ be a general model with signature $\Sigma$. Then for all $X \in \mathcal G$ and all variables $y$, $X[M/y] \in \mathcal G$.
\end{lemma}
\begin{proof}
	Let $\dom(X) = \{\tuple x\}$, let $R = \rel(X)$, and consider the formula $\phi(\tuple x, y) = \exists y R(\tuple x)$. Then take any assignment $s$ with domain $\{\tuple x, y\}$: by construction, $M \models_s \phi(\tuple x, y) \Leftrightarrow \exists m \mbox{ s.t. } s[m/y]_{|\tuple x} \in X \Leftrightarrow s \in X[M/y]$, as required.\footnote{Here by $s[m/y]_{|\tuple x}$ we intend the restriction of $s[m/y]$ to the domain $\{x_1 \ldots x_n\}$. If $y$ is among $x_1 \ldots x_n$, then this is the same of $s[m/y]$ itself; otherwise, it is simply $s$.}
\end{proof}

We can easily adapt the team semantics of the previous section to general models. We report all the rules here, for ease of reference; but the only differences between this semantics and the previous one are in the cases \textbf{PTS-}$\vee$ and \textbf{PTS-}$\exists$.
\begin{definition}[General Team Semantics for Independence Logic]
	Let $\Sigma$ be a first order signature, let $(M, \mathcal G)$ be a general model of signature $\Sigma$, let $\phi \in \nnf_\Sigma$ be a formula of Independence Logic and let $X \in \mathcal G$ be a team with domain containing $\free(\phi)$. Then we say that $X$ \emph{satisfies} $\phi$ in $(M, \mathcal G)$, and we write $(M, \mathcal G) \models_X \phi$, if and only if 
	\begin{description}
		\item[GTS-lit] $\phi$ is a first order literal and, for all $s \in X$, $M \models_s \phi$ in the usual first order sense; 
		\item[GTS-ind] $\phi$ is $\indep{\tuple t_1}{\tuple t_2}{\tuple t_3}$ for some tuples of terms $\tuple t_1$, $\tuple t_2$ and $\tuple t_3$, and for all $s, s' \in X$ with $\tuple t_1\langle s\rangle = \tuple t_1\langle s'\rangle$ there exists a $s'' \in X$ with $\tuple t_1 \tuple t_2\langle s''\rangle = \tuple t_1 \tuple t_2\langle s\rangle$ and $\tuple t_1 \tuple t_3\langle s''\rangle = \tuple t_1 \tuple t_3\langle s'\rangle$;
		\item[GTS-$\vee$] $\phi$ is $\psi_1 \vee \psi_2$ for two formulas $\psi_1, \psi_2 \in \nnf_{\Sigma}$ and $X = Y \cup Z$ for some two teams $Y, Z \in \mathcal G$ such that $(M, \mathcal G) \models_Y \psi_1$ and $(M, \mathcal G) \models_Z \psi_2$;
		\item[GTS-$\wedge$] $\phi$ is $\psi_1 \wedge \psi_2$ for two formulas $\psi_1, \psi_2 \in \nnf_{\Sigma}$, $(M, \mathcal G) \models_X \psi_1$ and $(M, \mathcal G) \models_X \psi_2$;
		\item[GTS-$\exists$] $\phi$ is $\exists x \psi$ for some variable $x$ and some $\psi \in \nnf_\Sigma$ and there exists a team $X' \in \mathcal G$ such that $X[x]X'$ and such that $(M, \mathcal G) \models_{X'} \psi$; 
		\item[GTS-$\forall$] $\phi$ is $\forall x \psi$ for some suitable $x$ and $(M, \mathcal G) \models_{X[M/x]} \psi$.
	\end{description}
\end{definition}
The usual semantics for Independence Logic satisfies a \emph{locality principle}: in brief, the satisfiability of a formula $\phi$ in a team depends only on the restriction of the team to $\free(\phi)$. Let us verify that the same holds for entailment semantics: 
\begin{lemma}
	Let $(M, \mathcal G)$ be a general model, and let $X \in \mathcal G$ be such that $\dom(X) = \tuple x \tuple y$. Then $X_{|\tuple x} = \{s : \dom(s) = \tuple x, \exists \tuple m \mbox{ s.t. } s[\tuple m/\tuple y] \in X\}$ is in $\mathcal G$. \\

	Furthermore, let $Y \subseteq X_{|\mathcal G}$ be such that $Y \in \mathcal G$. Then the team 
	\[
		X(\tuple x \in Y) = \{s \in X: s_{|\tuple x} \in Y\}
	\]
 	is in $\mathcal G$.
\end{lemma}
\begin{proof}
	By definition, $X_{|\tuple x}$ is $\|\phi(\tuple x, R)\|_M$, where $\phi$ is $\exists \tuple y (R\tuple x \tuple y)$ and $R = \rel(X)$. Therefore, $X_{|\tuple x} \in \mathcal G$.\\

	Similarly, $X(\tuple x \in Y)$ is $\|\phi(\tuple x \tuple y, R_1, R_2)\|_M$, where $\phi$ is $R_1 \tuple x \tuple y \wedge R_2 \tuple x$, $R_1$ is $\rel(X)$ and $R_2$ is $\rel(Y)$. 
\end{proof}
\begin{theorem}[Locality]
	\label{theo:local}
	Let $(M, \mathcal G)$ be a general model, let $X \in \mathcal G$ and let $\phi$ be an independence logic formula over the signature of $M$ with $\free(\phi) = \tuple z \subseteq \dom(X)$. Then $(M, \mathcal G) \models_X \phi$ if and only if $(M, \mathcal G) \models_{X_{|\tuple z}} \phi$. 
\end{theorem}
\begin{proof}
	The proof is by structural induction on $\phi$. We present only the passages corresponding to disjunction and existential quantification, as the others are trivial:
	\begin{itemize}
		\item Suppose that $(M, \mathcal G) \models_X \psi_1 \vee \psi_2$. Then, by definition, there exist teams $Y$ and $Z$ in $\mathcal G$ such that $X = Y \cup Z$, $(M, \mathcal G) \models_{Y} \psi_1$ and $M \models_Z \psi_2$. By induction hypothesis, this means that $(M, \mathcal G) \models_{Y_{|\tuple z}} \psi_1$ and $(M, \mathcal G) \models_{Z_{|\tuple z}} \psi_2$. But $Y_{|\tuple z} \cup Z_{|\tuple z} = X_{|\tuple z}$, and hence $(M, \mathcal G) \models_{X_{|\tuple z}} \psi_1 \vee \psi_2$. 

			Conversely, suppose that $(M, \mathcal G) \models_{X_{|\tuple z}} \psi_1 \vee \psi_2$. Then there exist teams $Y', Z'$ in $\mathcal G$ such that $(M, \mathcal G) \models_{Y'} \psi_1$, $(M, \mathcal G) \models_{Z'} \psi_2$ and $X_{|\tuple z} = X' \cup Y'$. Now let $Y$ be $X(\tuple z \in Y')$ and $Z$ be $X(\tuple z \in Z')$; by construction, $Y \cup Z = X$, and furthermore $Y' = Y_{|\tuple z}$ and $Z' = Z_{|\tuple z}$, and, by the lemma, $Y$ and $Z$ are in $\mathcal G$. 
			Thus, by induction hypothesis, $(M, \mathcal G) \models_{Y} \psi_1$ and $(M, \mathcal G) \models_Z \psi_2$, and finally $(M, \mathcal G) \models_X \psi_1 \vee \psi_2$, as required. 
		\item Suppose that $(M, \mathcal G) \models_X \exists x \psi$. Then there exists a team $Y \in \mathcal G$ such that $X[x]Y$ and $(M, \mathcal G) \models_{Y} \psi$. By induction hypothesis, this means that $(M, \mathcal G) \models_{Y_{|\tuple z x}} \psi$ too; and since $X_{|\tuple z}[x] Y_{|\tuple z x}$, this implies that $M \models_{X_{|\tuple z}} \exists x \psi$, as required. 

			Conversely, suppose that $(M, \mathcal G) \models_{X_{|\tuple z}} \exists x \psi$. Then there exists a team $Y'$, with domain $\tuple z x$, such that $M \models_{Y'} \psi$ and $X_{|\tuple z}[x]Y'$. Now let $Y$ be $(X[M/x])(\tuple z x \in Y')$. By the lemma, $Y \in \mathcal G$; furthermore, $Y_{|\tuple z x} = Y'$, and hence by induction hypothesis $(M, \mathcal G) \models_Y \psi$. Finally, $X[x]Y$: indeed, if $s \in X$ then $s_{\tuple z}[m/x] \in Y'$ for some $m \in \dom(M)$, and hence $s[m/x] \in Y$ for the same $m$, and on the other hand, $Y$ is contained in $X[M/x]$, and hence if $s[m/x] \in Y$ it follows that $s \in X$. 

			Therefore $(M, \mathcal G) \models_X \exists x \psi$, as required. 
	\end{itemize}
\end{proof}

As in the case of Second Order Logic, first-order models can be represented as a special kind of general model:
\begin{definition}[Full models]
	Let $(M, \mathcal G)$ be a general model. Then it is said to be \emph{full} if and only if $\mathcal G$ contains all teams over $M$. 
\end{definition}
The following result is then trivial. 
\begin{proposition}
	Let $(M, \mathcal G)$ be a full model. Then for all suitable teams $X$ and formulas $\phi$, $(M, \mathcal G) \models_X \phi$ in general team semantics if and only if $M \models_X \phi$ in the usual team semantics.
\end{proposition}
\begin{proof}
	Follows at once by comparing the rules of Team Semantics and General Team Semantics for the case that $\mathcal G$ contains all teams. 
\end{proof}

How does the satisfaction relation in general team semantics change if we vary the set $\mathcal G$? The following definition and result give us some information about this: 
\begin{definition}[Refinement]
	Let $(M, \mathcal G)$ and $(M, \mathcal G')$ be two general models. Then we say that $(M, \mathcal G')$ is a \emph{refinement} of $(M, \mathcal G)$, and we write $(M, \mathcal G) \subseteq (M, \mathcal G')$, if and only if $\mathcal G \subseteq \mathcal G'$. 
\end{definition}
Intuitively speaking, a refinement of a general model is another general model with more teams than it. The following result shows that refinements preserve satisfaction relations:
\begin{theorem}
	\label{theo:refin}
	Let $(M, \mathcal G)$ and $(M, \mathcal G')$ be two general models with\\ $(M, \mathcal G) \subseteq (M, \mathcal G')$, let $X \in \mathcal G$, and let $\phi$ be a formula over the signature of $M$ with $\free(\phi) \subseteq \dom(X)$. Then
	\[
		(M, \mathcal G) \models_X \phi \Rightarrow (M, \mathcal G') \models_X \phi.
	\]
\end{theorem}
\begin{proof}
	The proof is an easy induction on $\phi$. 
	\begin{enumerate}
		\item If $\phi$ is a first order literal, the result is obvious, as the choice of the set of teams $\mathcal G$ (or $\mathcal G'$) does not enter into the definition of satisfaction condition \textbf{PTS-lit}. 
		\item If $\phi$ is an independence atom, the result is also obvious, for the same reason.
		\item If $(M, \mathcal G) \models_X \psi_1 \vee \psi_2$ then there exist two teams $Y, Z \in \mathcal G$ such that $X = Y \cup Z$, $(M, \mathcal G) \models_Y \psi_1$ and $(M, \mathcal G) \models_Z \psi_2$. But $Y$ and $Z$ are also in $\mathcal G'$, and by induction hypothesis we have that $(M, \mathcal G') \models_Y \psi_1$ and $(M, \mathcal G') \models_Z \psi_2$, and therefore $(M, \mathcal G') \models_X \psi_1 \vee \psi_2$. 
		\item If $(M, \mathcal G) \models_X \psi_1 \wedge \psi_2$ then $(M, \mathcal G) \models_X \psi_1$ and $(M, \mathcal G) \models_X \psi_2$. Then, by induction hypothesis, $(M, \mathcal G') \models_X \psi_1$ and $(M, \mathcal G') \models_X \psi_2$, and finally $(M, \mathcal G') \models_X \psi_1 \wedge \psi_2$. 
		\item If $(M, \mathcal G) \models_X \exists x \psi$ then there exists a $X' \in \mathcal G$ such that $X[x]X'$ and $(M, \mathcal G) \models_{X'} \psi$. But then $X'$ is also in $\mathcal G'$, and by induction hypothesis $(M, \mathcal G') \models_{X'} \psi$, and finally $(M, \mathcal G') \models_X \exists x \psi$. 
		\item If $(M, \mathcal G) \models_X \forall x \psi$ then $(M, \mathcal G) \models_{X[M/x]} \psi$. Then, by induction hypothesis, $(M, \mathcal G') \models_{X[M/x]} \psi$, and finally $(M, \mathcal G') \models_X \forall x \psi$.
	\end{enumerate}
\end{proof}
This result shows us that, as was to be expected from the equivalence between independence logic and existential second order logic, if we are interested in formulas which hold in \emph{all} general models over a certain first-order model we only need to pay attention to the \emph{smallest} (in the sense of the refinement relation) ones. But do such ``least general models'' exist? As the following result shows, this is indeed the case: 
\begin{proposition}
	Let $\{(M, \mathcal G_i) : i \in I\}$ be a family of general models with signature $\Sigma$ and over the same first order model $M$. Then $(M, \bigcap_{i \in I} \mathcal G_i)$ is also a general model. 
\end{proposition}
\begin{proof}
		Let $\phi(x_1 \ldots x_n, \tuple m, \tuple R)$ be a first order formula with parameters, where each $R_i$ is of the form $\rel(X)$ for some $X \in \cap_i \mathcal G_i$. 
	Then the team $\|\phi(x_1 \ldots x_n, \tuple m, \tuple R)\|_M$ is in $\mathcal G_i$ for all $i \in I$, and therefore it is in $\bigcap_{i \in I} \mathcal G$, as required.
\end{proof}

Therefore, it is indeed possible to talk about the \emph{least general model} over a first order model. 
\begin{definition}[Least General Model]
	Let $M$ be a first order model. Then the \emph{least general model} over $M$ is the $(M, \mathcal L)$, where 
	\[
	\mathcal L = \bigcap \{\mathcal G : (M, \mathcal G) \mbox{ is a general model.}\}
	\]
\end{definition}
What is the purpose of least general models? The answer comes as a consequence of Theorem \ref{theo:refin}, and can be summarized by the following corollary: 
\begin{corollary}
	Let $\Sigma$ be a first order signature, let $M$ be a first order model over it and let $(M, \mathcal L)$ be the least general model over it. Then, for all teams $X \in \mathcal L$ and all formulas $\phi$ with signature $\Sigma$ and with free variables in $\dom(X)$, 
	\[
		(M, \mathcal L) \models_X \phi \Leftrightarrow (M, \mathcal G) \models_X \phi \mbox{ for all general models } (M, \mathcal G) \mbox{ over M}.
	\]
\end{corollary}
\begin{proof}
	Suppose that $(M, \mathcal L) \models_X \phi$. Then take any general model $(M, \mathcal G)$: by definition, we have that $(M, \mathcal L) \subseteq (M, \mathcal G)$, and hence by Theorem \ref{theo:refin} we have that $(M, \mathcal G) \models_X \phi$.

	Conversely, suppose that $(M, \mathcal G) \models_X \phi$ for all general models $(M, \mathcal G)$; then in particular $(M, \mathcal L) \models_X \phi$, as required. 
\end{proof}

We can also find a more practical characterization of this ``least general model''.
\begin{proposition}
	Let $M$ be a first order model. Then the least general model over it is $(M, \mathcal L)$, where $\mathcal L$ is the set of all $\|\phi(\tuple x, \tuple m)\|_M$, where $\phi$ ranges over all first order formulas and $\tuple m$ ranges over all tuples of variables of suitable length. 
\end{proposition}
\begin{proof}
	If $(M, \mathcal G)$ is a general model then $\mathcal L \subseteq \mathcal G$ by definition; therefore, we only need to prove that $(M, \mathcal L)$ is a general model.

	Now, let $\phi(\tuple x, \tuple m, \tuple R)$ be a first order formula, and let each $R_i$ be $\rel(X_i)$ for some $X_i \in \mathcal L$. So for each $R_i$, any assignment $s$ and any suitable tuple of terms $t$, $M \models_s R_i \tuple t$ if and only if $M \models_s \psi_i(\tuple t, \tuple n_i)$ for some first order formula $\psi_i$ with parameters $\tuple n_i$. Now let $\phi'(\tuple x, \tuple m, \tuple n_1, \tuple n_2, \ldots)$ be the expression obtained by substituting, in $\phi$, each instance of $R_i \tuple t$ with $\psi_i(\tuple t, \tuple n_i)$; by construction, we have that $M \models_s \phi(\tuple x, \tuple m, \tuple R)$ if and only if $M \models_s \phi'(\tuple x, \tuple m, \tuple n_1, \ldots)$, and therefore 
	\[
		\|\phi(\tuple x, \tuple m, \tuple R)\|_M = \|\phi'(\tuple x, \tuple m, \tuple n_1, \tuple n_2, \ldots)\|_M \in \mathcal L
	\]
	as required. 
\end{proof}

As long as we are only considering teams in $\mathcal L$, studying satisfiability with respect to the least model $(M, \mathcal L)$ is the same as considering satisfiability with respect to all general models over $M$. This restriction may at first sight seem a bit unpractical, but it becomes irrelevant when it comes to the problem of \emph{validity}:
\begin{definition}[Validity wrt general models]
	Let $\Sigma$ be a first order signature, let $V$ be a finite set of variables, and let $\phi \in \nnf_\Sigma$ be a formula of our language with free variables in $V$. Then $\phi$ is \emph{valid} with respect to general models if and only if $(M, \mathcal G) \models_X \phi$ for all general models $(M, \mathcal G)$ with signature $\Sigma$ and for all teams $X \in \mathcal G$ with $\dom(X) \supseteq \free(\phi)$. If this is the case, we write $\gts \models \phi$. 
\end{definition}
\begin{definition}[Validity wrt least general models]
	Let $\Sigma$ be a first order signature, let $V$ be a finite set of variables, and let $\phi \in \nnf_\Sigma$ be a formula of our language with free variables in $V$. Then $\phi$ is \emph{valid} with respect to least general models if and only if $(M, \mathcal L) \models_X \phi$ for all least general models $(M, \mathcal L)$ with signature $\Sigma$ and for all teams $X \in \mathcal L$ with $\dom(X) \supseteq \free(\phi)$. If this is the case, we write $\lts \models \phi$. 
\end{definition}
\begin{lemma}
	\label{lemma:restrict_sig}
	Let $M$ be a first order model with signature $\Sigma$, and let $M'$ be another first order model with signature $\Sigma' \supseteq \Sigma$ such that the restriction of $M'$ to $\Sigma$ is precisely $M$. Then for all general models $\mathcal G$ for $M'$, for all formulas $\phi$ with signature $\Sigma$ and for all $X \in \mathcal G$, 
	\[
		(M, \mathcal G) \models_X \phi \Leftrightarrow (M', \mathcal G) \models_X \phi.
	\]
\end{lemma}
\begin{proof}
	First of all, if $(M', \mathcal G)$ is a general model then $(M, \mathcal G)$ is also a general model. Then, the result is proved by observing that the truth conditions of our semantics depend only on the interpretations of the symbols in the signature of the formula (and on the choice of $\mathcal G$, of course). 
\end{proof}
\begin{lemma}
	\label{lemma:enlarge_sig}
	Let $(M, \mathcal G)$ be a general model with signature $\Sigma$, let $S \not \in \Sigma$ be a new relation symbol and let $X \in \mathcal G$. Furthermore, let $M' = M[\rel(X)/S]$ be the extension of $M$ to the signature $\Sigma \cup \{S\}$ such that $S^{M'} = \rel(X)$. Then $(M', \mathcal G)$ is a general model. 
\end{lemma}
\begin{proof}
	Let $\phi(\tuple x, \tuple m, \tuple R)$ be a first order formula with signature $\Sigma \cup \{S\}$ and parameters $\tuple m$ and $\tuple R$, where each $R_i$ is $\rel(X_i)$ for some $X_i \in \mathcal G$. Then let $\phi'(\tuple x, \tuple m, \tuple R, S)$ be the first order formula with signature $\Sigma$, where $S$ now stands for the relation $\rel(X)$. Now clearly
	\[
		\|\phi(\tuple x, \tuple m, \tuple R)\|_{M'} = \|\phi'(\tuple x, \tuple m, \tuple R, S)\|_{M} \in \mathcal G,
	\]
	as required.
\end{proof}
\begin{theorem}
	\label{theo:valid_least}
	A formula $\phi$ is valid wrt general models if and only if it is valid wrt least general models. 
\end{theorem}
\begin{proof}
	The left to right direction is obvious. For the right to left direction, suppose that $\lts \models \phi$, let $(M, \mathcal G)$ be a general model whose signature contains the signature of $\phi$, and let $X \in \mathcal G$ be a team whose domain $\{x_1 \ldots x_n\}$ contains all free variables of $\phi$. Then consider the first order model $M' = M[\rel(X)/S]$, where $S$ is a new relation symbol, and take the least general model $(M', \mathcal L)$ over it. We clearly have that $X \in \mathcal L$, since 
	\[
	X = \{s: \dom(s) = \{x_1 \ldots x_n\}, M' \models_s S x_1 \ldots x_n\}
	\]
	and, therefore, $(M', \mathcal L) \models_X \phi$ by hypothesis. Now, by Lemma \ref{lemma:enlarge_sig}, $(M', \mathcal G)$ is a general model, and therefore by definition $\mathcal L \subseteq \mathcal G$, and hence by Theorem \ref{theo:refin} $(M', \mathcal G) \models_X \phi$ too. Finally, the relation symbol $S$ does not occur in $\phi$, and therefore by Lemma \ref{lemma:restrict_sig} $(M, \mathcal G) \models_X \phi$, as required.
\end{proof}

In the next section, we will develop another, more syntactic way of reasoning about least general models. 
\section{Entailment semantics}
Let $(M, \mathcal L)$ be a least general model; then, as we saw, $\mathcal L$ is the set of all teams corresponding to first order formulas with parameters. Therefore, in order to reason about satisfaction in a least general team, there is no need to carry around sets of assignments; rather, we can use these first order formulas themselves as our context-carrying objects. In this section, we will develop this idea, building up a new ``entailment semantics'' and proving its equivalence with ``least general model semantics''. \\

In order to do all of this, we first need to be a bit more formal about the concept of ``constant parameters''. This is achieved through the following definitions: 
\begin{definition}[Parameter and Team Variables]
	Let $\mathbf V_P = \{p_1 \ldots p_n, \ldots\}$ be a fixed, countably infinite set of variables. We will call them \emph{parameter variables}. The variables $\mathbf V_T = \{x, y, z, \ldots\} = \mathbf V \backslash V_P$ will be called \emph{team variables}.
\end{definition}
\begin{definition}[Free Parameter and Team Variables]
	Let $\phi$ be any formula. Then $\free_P (\phi) = \free(\phi) \cap \mathbf V_P$ and $\free_T(\phi) = \free(\phi) \cap \mathbf V_T$.
\end{definition}
Parameter variables clarify the interpretation of such expressions such as $M \models_s \gamma(\tuple x, \tuple m)$: this is simply a shorthand $M \models_{h \cup s} \gamma(\tuple x, \tuple p)$, where $h$ is a \emph{parameter assignment} with domain $\tuple p$ and with $h(\tuple p) = \tuple m$. Team variables, instead, are going to be used in order to describe the variables in the domain of the team corresponding to a given first order expression: for any first order $\gamma(\tuple x, \tuple p)$, where $\tuple x$ are team variables and $\tuple p$ are parameter variables, and for any $h$ with domain $\tuple p$, we will therefore have $\|\gamma(\tuple x, \tuple p)\|_{M, h} = \|\gamma(\tuple x, h(\tuple p))\|_M = \{s : \dom(s) = \tuple x, M \models_{h \cup s} \gamma\}$. For this reason, parameter variables will never occur in the domain of a team, and, hence, from this point on we will always assume that parameter variables never occur in independence logic formulas, but only in the first order team definitions. \\

After these preliminaries, we can now give our main definition for this section:
\begin{definition}[Entailment Semantics for Independence Logic]
	Let $M$ be a first order model with signature $\Sigma$, let $\gamma(\tuple x, \tuple p)$ be a first order formula for the same signature with $\free_T = \tuple x$ and $\free_P = \tuple p$, let $h$ be a parameter assignment with domain $\tuple p$ and let $\phi \in \nnf_{\Sigma}$ be an Independence Logic formula. \\
	
	Then we say that $\gamma$ \emph{satisfies} $\phi$ in $M$ \emph{under} $h$, and we write $M \models_{\gamma(h)} \phi$, if and only if 
	\begin{description}
		\item[ES-lit] $\phi$ is a first order literal and for all assignments $s$ with domain $\free_T(\gamma) \cup \free_T(\phi)$ such that $M \models_{h \cup s} \gamma$ it holds that $M \models_s \phi$;
		\item[ES-ind] $\phi$ is $\indep{\tuple t_1}{\tuple t_2}{\tuple t_3}$ for some tuples of terms $\tuple t_1$, $\tuple t_2$ and $\tuple t_3$ and for all assignments $s$ and $s'$ with domain $\free_T(\gamma) \cup \free_T(\tuple t_1 \tuple t_2 \tuple t_3)$ such that $M \models_{h \cup s} \gamma$, $M \models_{h \cup s'} \gamma$ and $\tuple t_1 \langle s\rangle = \tuple t_1 \langle s'\rangle$ there exists a $s''$ such that $M \models_{h \cup s''} \gamma$, $\tuple t_1 \tuple t_2\langle s''\rangle = \tuple t_1 \tuple t_2\langle s\rangle$ and $\tuple t_1 \tuple t_3\langle s''\rangle = \tuple t_1\tuple t_3\langle s'\rangle$;
		\item[ES-$\vee$] $\phi$ is $\psi_1 \vee \psi_2$ and there exists a parameter assignment $h'$ extending\footnote{That is, $\dom(h') \supseteq \dom(h)$, and $h'(\tuple p) = h(\tuple p)$.} $h$ and two first order formulas $\gamma_1$ and $\gamma_2$ such that 
			\begin{itemize}
				\item $\free_P(\gamma_1), \free_P(\gamma_2) \subseteq \dom(h')$;
				\item $M \models_{\gamma_1(h')} \psi_1$; 
				\item $M \models_{\gamma_2(h')} \psi_2$; 
				\item $M \models_{h'} \forall \tuple v(\gamma \leftrightarrow \gamma_1 \vee \gamma_2)$, where $\tuple v$ is $\free_T(\gamma) \cup \free_T(\gamma_1) \cup \free_T(\gamma_2)$; 
			\end{itemize}
		\item[ES-$\wedge$] $\phi$ is $\psi_1 \wedge \psi_2$, $M \models_{\gamma(h)} \psi_1$ and $M \models_{\gamma(h)} \psi_2$; 
		\item[ES-$\exists$] $\phi$ is $\exists x_n \psi$ and there exist a parameter assignment $h'$ extending $h$ and a first order formula $\gamma'$ with $\free_P(\gamma') \subseteq \dom(h')$ such that 
			\begin{itemize}
				\item $M \models_{\gamma'(h')} \psi$;
				\item $M \models_{h'} \forall \tuple v(\exists x_n \gamma' \leftrightarrow \exists x_n \gamma)$, where $\tuple v$ is $\free_T(\gamma) \cup \free_T(\gamma')$;
			\end{itemize}
		\item[ES-$\forall$] $\phi$ is $\forall x_n \psi$ and there exists a parameter assignment $h'$ extending $h$ and a first order formula $\gamma'$ with $\free_P(\gamma') \subseteq \dom(h')$ such that
						\begin{itemize}
				\item $M \models_{\gamma'(h')} \psi$;
				\item $M \models_{h'} \forall \tuple v(\gamma' \leftrightarrow \exists x_n \gamma)$, where $\tuple v$ is $\free_T(\gamma) \cup \free_T(\gamma')$.
			\end{itemize}
	\end{description}
\end{definition}
\begin{proposition}
	\label{propo:non_par_irrel}
	Let $M$ be a first order model with signature $\Sigma$, let $\gamma(\tuple x, \tuple p)$ be a first order formula with $\free_P(\gamma) =  \tuple p$ and let $h$, $h'$ be two parameter assignments with domains containing $\tuple p$ such that $h(\tuple p) = h'(\tuple p)$. Then, for all independence logic formulas $\phi$,
	\[
	M \models_{\gamma(h)} \phi \Leftrightarrow M \models_{\gamma(h')} \phi.
	\]
\end{proposition}
\begin{proof}
	The proof is a straightforward induction over $\phi$. 
\end{proof}

As the next result shows, entailment semantics is entirely equivalent to least general model semantics: 
\begin{theorem}
	\label{theo:ent_least}
	Let $\Sigma$ be a first order model, let $\gamma(\tuple x, \tuple p)$ be a first order formula with $\free_P(\gamma) = \tuple p$, let $h$ be a parameter assignment with domain $\tuple p$ and let $\phi \in \nnf_\Sigma$ be an independence logic formula with free variables in $\tuple x$.\\
	
	Furthermore, let $(M, \mathcal L)$ be the least general model over $M$, and let $X = \|\gamma(\tuple x, \tuple p)\|_{M, h} = \{s : \dom(s) = \{\tuple x\}, M \models_{h \cup s} \gamma(\tuple x, \tuple m)\}$. Then 
	\[
	(M, \mathcal L) \models_X \phi \Leftrightarrow M \models_{\gamma(h)} \phi.
	\]
\end{theorem}
\begin{proof}
	The proof is by structural induction on $\phi$, and presents no difficulties.
	\begin{enumerate}
		\item If $\phi$ is a first order literal, $(M, \mathcal L) \models_X \phi$ if and only if, for all $s \in X$, it holds that $M \models_s \phi$. But $s \in X$ if and only if $M \models_s \gamma(\tuple x, h(\tuple p))$, and hence $(M, \mathcal L) \models_X \phi$ if and only if $M \models_{\gamma} \phi$, as required.
		\item If $\phi$ is an independence atom, the result is also obvious, and follows at once from a comparison of the rules \textbf{PTS-ind} and \textbf{ES-ind}.
		\item If $\phi$ is $\psi_1 \vee \psi_2$,
			\begin{align*}
				&(M, \mathcal L) \models_{X} \psi_1 \vee \psi_2 \Leftrightarrow\\
				&\Leftrightarrow \exists Y, Z \in \mathcal L \mbox{ s.t. } X = Y \cup Z, (M, \mathcal L) \models_Y \psi_1 \mbox{ and } (M, \mathcal L) \models_Z \psi_2 \Leftrightarrow\\
				&\Leftrightarrow \exists h' = h[\tuple m/\tuple q] \mbox{ extending } h \mbox{ and } \exists \gamma_1 \gamma_2 \mbox{ s.t., for } Y = \|\gamma_1(\tuple x, \tuple p \tuple q)\|_{M, h'},\\
	&Z = \|\gamma_2(\tuple x, \tuple p \tuple q)\|_{M, h'},  X = \|\gamma(\tuple x, \tuple p)\|_{M, h} = \|\gamma(\tuple x, \tuple p)\|_{M, h'} = Y \cup Z,\\
	&(M, \mathcal L) \models_{Y} \psi_1 \mbox{ and }(M, \mathcal L) \models_{Z} \psi_2 \Leftrightarrow\\
				&\Leftrightarrow \exists h' = h[\tuple m/\tuple q] \mbox{ extending } h \mbox{ and } \exists \gamma_1 \gamma_2 \mbox{ s.t. } M \models_{h'} \forall \tuple v( \gamma \leftrightarrow \gamma_1 \vee \gamma_2),\\
	& M \models_{\gamma_1(h')} \psi \mbox{ and } M \models_{\gamma_2(h')} \theta \Leftrightarrow\\
				&\Leftrightarrow M \models_{\gamma(h)} \psi \vee \theta.
			\end{align*}
		\item If $\phi$ is $\psi \wedge \theta$,
			\begin{align*}
				&(M, \mathcal L) \models_{X} \psi \wedge \theta \Leftrightarrow (M, \mathcal L) \models_{X} \psi \mbox{ and } (M, \mathcal L) \models_{X} \theta \Leftrightarrow\\
				&\Leftrightarrow M \models_{\gamma(h)} \psi \mbox{ and } M \models_{\gamma(h)} \theta \Leftrightarrow M \models_{\gamma(h)} \psi \wedge \theta.
			\end{align*}
		\item If $\phi$ is $\exists x_n \psi$, 
			\begin{align*}
				&(M, \mathcal L) \models_{X} \exists x_n \psi \Leftrightarrow \exists X' \in \mathcal L \mbox{ s.t. } X[x_n]X' \mbox{ and } (M, \mathcal L) \models_{X'}\psi \Leftrightarrow\\
				&\Leftrightarrow \exists h' = h[\tuple m/\tuple q] \mbox{ extending } h \mbox{ and } \exists \gamma' \mbox{ s.t., for } X' = \|\gamma'(\tuple x, \tuple p \tuple q)\|_{M, h'},\\
& X[x_n]X' \mbox{ and } (M, \mathcal L) \models_{X'} \psi \Leftrightarrow\\
				&\Leftrightarrow \exists h' = h[\tuple m/\tuple q] \mbox{ extending } h \mbox{ and } \exists \gamma' \mbox{ s.t. } M \models_{h'} \forall \tuple v(\exists x_n \gamma \leftrightarrow \exists x_n \gamma') \mbox{ and}\\
& \mbox{ and } M \models_{\gamma'(h')} \psi \Leftrightarrow\\
				&\Leftrightarrow M \models_{\gamma(h)} \exists x_n \psi;
			\end{align*}
		\item If $\phi$ is $\forall x_n \psi$, 
				\begin{align*}
				&(M, \mathcal L) \models_{X} \forall x_n \psi \Leftrightarrow \exists X' \in \mathcal L \mbox{ s.t. } X'=X[M/x_n] \mbox{ and } (M, \mathcal L) \models_{X'}\psi \Leftrightarrow\\
				&\Leftrightarrow \exists h' = h[\tuple m/\tuple q] \mbox{ extending } h \mbox{ and } \exists \gamma' \mbox{ s.t., for } X' = \|\gamma'(\tuple x, \tuple pq)\|_{M, h'},\\
				& X' = X[M/x_n] \mbox{ and } (M, \mathcal L) \models_{X'} \psi \Leftrightarrow\\
				&\Leftrightarrow \exists h' = h[\tuple m/\tuple q] \mbox{ extending } h \mbox{ and } \exists \gamma' \mbox{ s.t. } M \models_{h'} \forall \tuple v(\gamma' \leftrightarrow \exists x_n \gamma) \mbox{ and}\\
	&\mbox{ and } M \models_{\gamma'(h')} \psi \Leftrightarrow\\
				&\Leftrightarrow M \models_{\gamma(h)} \forall x_n \psi.
			\end{align*}
	\end{enumerate}
\end{proof}
\begin{definition}[Validity in Entailment Semantics]
	Let $\phi$ be an Independence Logic formula. Then $\phi$ is \emph{valid} in entailment semantics if and only if $M \models_{\gamma(h)} \phi$ for all first order models $M$ with signature containing that of $\phi$, for all first order formulas $\gamma(\tuple x, \tuple p)$ over the signature of $M$ and for all parameter assignments $h$ with domain $\tuple p$. If this is the case, we write $\ens \models \phi$. 
\end{definition}
\begin{corollary}
	For all formulas $\phi$,  $\ens \models \phi$ if and only if $\lts \models \phi$ if and only if $\gts \models \phi$
\end{corollary}

It will also be useful to have a slightly more general notion of validity in entailment semantics: 
\begin{definition}[Validity wrt a Team Definition]
		Let $\gamma(\tuple x, \tuple p)$ be a first order formula and let $\phi$ be an independence logic formula. Then $\phi$ is \emph{valid} with respect to $\gamma$ if and only if $M \models_{\gamma(h)} \phi$ for all first order models $M$ with signature containing those of $\gamma$ and $\phi$ and for all parameter assignments $h$ with domain $\tuple p$. If this is the case, we write $\models_\gamma \phi$. 
\end{definition}
\begin{proposition}
	Let $\phi$ be an independence logic formula with $\free_T(\phi) = \tuple x$, and let $R$ be a $|\tuple x|$-ary relation symbol not occurring in $\gamma$. Then $\ens \models \phi$ if and only if $\models_{R\tuple x} \phi$.
\end{proposition}
\begin{proof}
	Suppose that $\ens \models \phi$. Then in particular, for any model $M$ whose signature contains that of $\phi$ and $R$ we have that $M \models_{R \tuple x} \phi$, and hence $\models_{R \tuple x} \phi$. \\

	Conversely, suppose that $\models_{R \tuple x} \phi$, let $M$ be a first order model\footnote{Without loss of generality, we can assume that the signature of $M$ does not contain the symbol $R$.}, and let $X \in \mathcal L$ be any team with domain $\tuple x$. Let us then consider the model $M'$ obtained by adding to $M$ the $|\tuple x|$-ary symbol $R$ with $R^{M'} = \rel(X)$. By hypothesis, $M' \models_{R \tuple x} \phi$, and furthermore since $R^{M'}$ is in $\mathcal L$ already the least general model over $M'$ is $(M', \mathcal L)$ for the same $\mathcal L$.
	
	Now $(M', \mathcal L) \models_{X} \phi$, and therefore, as $R$ occurs nowhere in $\phi$, $(M, \mathcal L) \models_X \phi$ too. This holds for all $X$ with domains $\tuple x$; therefore by the Locality Theorem (Theorem \ref{theo:local}), the same holds for all domains containing $\tuple x$, and hence $\lts \models \phi$. This implies that $\ens \models \phi$, as required.
\end{proof}

In the next section, we will develop a sound and complete proof system for this notion of validity with respect to a team definition.\\
\section{The proof system}
In this section, we will develop a proof system for Independence Logic (with entailment semantics) and prove its soundness and completeness. 

\begin{definition}[Sequent]
	Let $\Gamma$ be a finite first order theory with only parameter variables among its free ones, let $\gamma(\tuple x, \tuple p)$ be a first order formula and let $\phi$ be an Independence Logic formula with free variables in $\mathbf V_T$. Then the expression 
	\[
		\Gamma ~|~ \gamma \vdash \phi
	\]
	is a sequent.
\end{definition}
The intended semantics of a sequent is the following one: 
\begin{definition}[Valid Sequents]
	Let $\Gamma ~|~ \gamma \vdash \phi$ be a sequent. Then $\Gamma ~|~ \gamma \vdash \phi$ is \emph{valid} if and only if for all models $M$ and all parameter assignments $h$ with domain $\free_P(\Gamma) \cup \free_P(\gamma)$ such that $M \models_h \Gamma$ it holds that 
	\[
		M \models_{\gamma(h)} \phi.
	\]
\end{definition}
The following result is then clear: 
\begin{proposition}
	For all $\gamma$ and $\phi$, $\models_\gamma \phi$ if and only if $\emptyset ~|~ \gamma \vdash \phi$ is valid.
\end{proposition}

Now, all we need to do is develop some syntactic rules for finding whether a given sequent is valid. 

We can do this as follows:
\begin{definition}[Axioms and Rules]
	The axioms of our proof system are
	\begin{description}
		\item[PS-lit] If $\phi$ is a first order literal with no free parameter variables (that is, $\free_P(\phi) = \emptyset$) then
			\[
				\forall \tuple v (\gamma \rightarrow \phi) ~|~ \gamma \vdash \phi
			\]
			for all first order formulas $\gamma$, where $\tuple v = \free_T(\gamma) \cup \free_T(\phi)$;
		\item[PS-ind] If $\tuple t_1$, $\tuple t_2$ and $\tuple t_3$ are first order terms with no free parameter variables then 
			\begin{align*}
			&\forall \tuple v_1 \tuple v_2 ((\gamma(v_1) \wedge \gamma(v_2) \wedge \tuple t_1(\tuple v_1) = \tuple t_1(\tuple v_2)) \rightarrow \exists \tuple v_3(\gamma(v_3) \wedge \tuple t_1 \tuple t_2(\tuple v_3) = \tuple t_1 \tuple t_2(\tuple v_1) \wedge\\
			& ~ ~ ~ ~ ~ ~ \tuple t_1 \tuple t_3(\tuple v_3) = \tuple t_1 \tuple t_3(\tuple v_2))) ~|~ \gamma \vdash \indep{\tuple t_1}{\tuple t_2}{\tuple t_3}
			\end{align*}
			for all $\gamma$, where $\tuple v_1$ and $\tuple v_2$ are tuples of variables of the same lengths of $\tuple v = \free_T(\gamma) \cup \free_T(\tuple t_1 \tuple t_2 \tuple t_3)$, $\tuple t_i(\tuple v_j)$ is the tuple obtained by substituting $\tuple v$ with $\tuple v_j$ in $\tuple t_i$, and the same holds for $\gamma(\tuple v_j)$.
	\end{description}
	The rules of our proof system are
	\begin{description}
		\item[PS-$\vee$] If $\Gamma_1 ~|~ \gamma_1 \vdash \phi_1$ and $\Gamma_2 ~|~ \gamma_2 \vdash \phi_2$ then, for all $\gamma$, we have
			\[
				\Gamma_1, \Gamma_2, \forall \tuple v(\gamma \leftrightarrow (\gamma_1 \vee \gamma_2)) ~|~ \gamma \vdash \phi_1 \vee \phi_2
			\]
			where $\tuple v$ is $\free_T(\gamma) \cup \free_T(\gamma_1) \cup \free_T(\gamma_2)$; 
		\item[PS-$\wedge$] If $\Gamma_1 ~|~ \gamma \vdash \phi_1$ and $\Gamma_2 ~|~ \gamma \vdash \phi_2$ then $\Gamma_1, \Gamma_2 ~|~ \gamma \vdash \phi_1 \wedge \phi_2$; 
		\item[PS-$\exists$] If $\Gamma ~|~ \gamma' \vdash \phi$ and $x$ is a team variable then, for all $\gamma$,
			\[
				\Gamma, \forall \tuple v (\exists x \gamma' \leftrightarrow \exists x \gamma) ~|~ \gamma \vdash \exists x \phi
			\]
			where $\tuple v = \free_T(\gamma) \cup \free_T(\gamma')$; 
		\item[PS-$\forall$] If $\Gamma ~|~ \gamma' \vdash \phi$ and $x$ is a team variable then, for all $\gamma$, 
			\[
				\Gamma, \forall \tuple v (\gamma' \leftrightarrow \exists x \gamma) ~|~ \gamma \vdash \forall x \phi
			\]
			where, as in the previous case, $\tuple v = \free_T(\gamma) \cup \free_T(\gamma')$; 
		\item[PS-ent] If $\Gamma ~|~ \gamma \vdash \phi$ and $\bigwedge \Gamma' \models \bigwedge \Gamma$ holds in First Order Logic then $\Gamma' ~|~ \gamma \vdash \phi$;
		\item[PS-depar] If $\Gamma ~|~ \gamma \vdash \phi$ and $p$ is a parameter variable which does not occur free in $\gamma$ then $\exists p \bigwedge \Gamma ~|~ \gamma \vdash \phi$;
		\item[PS-split] If $\Gamma_1 ~|~ \gamma \vdash \phi$ and $\Gamma_2 ~|~ \gamma \vdash \phi$ then $(\bigwedge \Gamma_1) \vee (\bigwedge \Gamma_2) ~|~ \gamma \vdash \phi$.	
	\end{description}
\end{definition}
\begin{definition}[Proofs and proof lengths]
	Let $\Gamma ~|~ \gamma \vdash \phi$ be a sequent. A \emph{proof} of this sequent is a finite list of sequents 
	\[
		(\Gamma_1 ~|~ \gamma_1 \vdash \phi_1), \ldots, (\Gamma_n ~|~ \gamma_n \vdash \phi_n) = (\Gamma ~|~ \gamma \vdash \phi)
	\]
	such that, for all $i = 1 \ldots n$, $\Gamma_i ~|~ \gamma_i \vdash \phi_i$ is either an instance of \textbf{PS-lit}, \textbf{PS-ind} or it follows from $\{\Gamma_j ~|~ \gamma_j \vdash \phi_j : j < i\}$ through one application of the rules of our proof system. 

	Given a proof $P = S_1 \ldots S_n$, where each $S_i$ is a sequent, we define its \emph{length} $|P|$ as $n - 1$, that is, as the number of sequents in the proof minus one.
\end{definition}
Before examining soundness and completeness for this proof system, it will be useful to obtain a couple of derived rules:
\begin{proposition}
	The following rules hold:
	\begin{description}
		\item[PS-FO] If $\phi$ is a first order formula with no free parameter variables $\forall \tuple v(\gamma \rightarrow \phi) ~|~ \gamma \vdash \phi$  is provable for all $\gamma$, where $\tuple v = \free_T(\gamma) \cup \free_T(\phi)$; 
		\item[PS-dep] If $\tuple t$ is a tuple of terms, $t'$ is another term and $=\!\!(\tuple t, t')$ stands for $\indep{\tuple t}{t'}{t'}$ then 
			\[
				\forall \tuple v_1 \tuple v_2 (\gamma(\tuple v_1) \wedge \gamma( \tuple v_2) \wedge \tuple t(\tuple v_1) = \tuple t(\tuple v_2)) \rightarrow t'(\tuple v_1) = t'(\tuple v_2) ~|~ \gamma \vdash =\!\!(\tuple t, t')
			\]
			is provable for all $\gamma$, where $\tuple v_1$, $\tuple v_2$ are tuples of the same length of $\tuple v = \free_T(\gamma \cup \free_T(\tuple t t')$.
	\end{description}
\end{proposition}
\begin{proof}
	\begin{description}
		\item[PS-FO] The proof is by structural induction on $\phi$. 
			\begin{enumerate}
				\item If $\phi$ is a first order literal, this follows at once from rule \textbf{PS-lit}. 
				\item If $\phi$ is $\psi_1 \vee \psi_2$, by induction hypothesis we have that\\
 $\forall \tuple v ( (\gamma \wedge \psi_1) \rightarrow \psi_1) ~|~ \gamma \wedge \psi_1 \vdash \psi_1$ and $\forall \tuple v ((\gamma \wedge \psi_2) \rightarrow \psi_2) ~|~ \gamma \wedge \psi_2 \vdash \psi_2$ are provable. But then we can prove $\forall \tuple v(\gamma \rightarrow \phi_1 \vee \phi_2) ~|~ \gamma \vdash \phi$ as follows: 
					\begin{enumerate}
					\item $\forall \tuple v ( (\gamma \wedge \psi_1) \rightarrow \psi_1) ~|~ \gamma \wedge \psi_1 \vdash \psi_1$ (Derived before)
					\item $\forall \tuple v ((\gamma \wedge \psi_2) \rightarrow \psi_2) ~|~ \gamma \wedge \psi_2 \vdash \psi_2$ (Derived before)
					\item $~|~ \gamma \wedge \psi_1 \vdash \psi_1$ (\textbf{PS-ent}, from (a), because $\models \forall \tuple v ( (\gamma \wedge \psi_1) \rightarrow \psi_1)$ in First Order Logic)
					\item $~|~ \gamma \wedge \psi_2 \vdash \psi_2$ (\textbf{PS-ent}, from (b), because $\models \forall \tuple v ( (\gamma \wedge \psi_2) \rightarrow \psi_2)$ in First Order Logic)
					\item $\forall \tuple v(\gamma \leftrightarrow (\gamma \wedge \psi_1) \vee (\gamma \wedge \psi_2)) ~|~ \gamma \vdash \psi_1 \vee \psi_2$ (\textbf{PS-}$\vee$, from (c) and (d))
					\item $\forall \tuple v(\gamma \rightarrow (\psi_1 \vee \psi_2)) ~|~ \gamma \vdash \psi_1 \vee \psi_2$ (\textbf{PS-ent}: from (e), because $\forall \tuple v(\gamma \rightarrow (\psi_1 \vee \psi_2))$ entails $\forall \tuple v(\gamma \leftrightarrow (\gamma \wedge \psi_1) \vee (\gamma \wedge \psi_2))$ in First Order Logic).
					\end{enumerate}
				%
				\item If $\phi$ is $\psi_1 \wedge \psi_2$, by induction hypothesis we have that $\forall \tuple v(\gamma \rightarrow \psi_1) ~|~ \gamma \vdash \psi_1$ and $\forall \tuple v(\gamma \rightarrow \psi_2) ~|~ \gamma \vdash \psi_2$ are provable. But then
					\begin{enumerate}
						\item $\forall \tuple v(\gamma \rightarrow \psi_1) ~|~ \gamma \vdash \psi_1$ (derived before)
						\item $\forall \tuple v(\gamma \rightarrow \psi_2) ~|~ \gamma \vdash \psi_2$ (derived before)
						\item $\forall \tuple v(\gamma \rightarrow \psi_1), \forall \tuple v(\gamma \rightarrow \psi_2) ~|~ \gamma \vdash \psi_1 \wedge \psi_2$ (\textbf{PS}-$\wedge$, (a), (b))
						\item $\forall \tuple v(\gamma \rightarrow \psi_1 \wedge \psi_2) ~|~ \gamma \vdash \psi_1 \wedge \psi_2$ (\textbf{PS-ent}, (c))
					\end{enumerate}
					as required.	
				\item If $\phi$ is $\exists x \psi$, by induction hypothesis we have that\\ 
$\forall \tuple v \forall x ((\exists x \gamma) \wedge \psi) \rightarrow \psi) ~|~ (\exists x \gamma) \wedge \psi \vdash \psi$ is provable. But then 
					\begin{enumerate}
						\item $\forall \tuple v \forall x (((\exists x \gamma) \wedge \psi) \rightarrow \psi) ~|~ (\exists x \gamma) \wedge \psi \vdash \psi$ (derived before)
						\item $~|~ (\exists x \gamma) \wedge \psi \vdash \psi$ (\textbf{PS-ent}, from (a))
						\item $\forall \tuple v(\exists x ( (\exists x \gamma) \wedge \psi) \leftrightarrow \exists x \gamma) ~|~ \gamma \vdash \exists x \psi$ (\textbf{PS-}$\exists$, from (b))
						\item $\forall \tuple v(((\exists x \gamma) \wedge (\exists x \psi)) \leftrightarrow \exists x \gamma) ~|~ \gamma \vdash \exists x \psi$ (\textbf{PS-ent}, from (c))
						\item $\forall \tuple v (\gamma \rightarrow \exists x \psi) ~|~ \gamma \vdash \psi$ (\textbf{PS-ent}, from (d))
					\end{enumerate}
					as required, where the last passage uses the fact that\\ $\forall \tuple v (\gamma \rightarrow \exists x \psi) \models \forall \tuple v (((\exists x \gamma) \wedge (\exists x \psi)) \leftrightarrow \exists x \gamma)$ in First Order Logic.
				\item If $\phi$ is $\forall x \psi$, by induction hypothesis we have that\\ $\forall \tuple v \forall x ((\exists x \gamma) \rightarrow \psi) ~|~ \exists x \gamma \vdash \psi$ is provable. But then 
					\begin{enumerate}
						\item $\forall \tuple v\forall x ((\exists x \gamma) \rightarrow \psi) ~|~ \exists x \gamma \vdash \psi$ (derived before)
						\item $\forall \tuple v\forall x ((\exists x \gamma) \rightarrow \psi), \forall \tuple v(\exists x \gamma \leftrightarrow \exists x \gamma) ~|~ \gamma \vdash \forall x \psi$ (\textbf{PS-}$\forall$, from (a))
						\item $\forall \tuple v \forall x( (\exists x \gamma) \rightarrow \psi) ~|~ \gamma \vdash \forall x \psi$ (\textbf{PS-ent}, from (c))
						\item $\forall \tuple v (\gamma \rightarrow \forall x \psi) ~|~ \gamma \vdash \forall x \psi$ (\textbf{PS-ent}, from (d))
					\end{enumerate}
					where the last two passages hold because $\forall \tuple v(\exists x \gamma \leftrightarrow \exists x \gamma)$ is valid and because $\forall \tuple v (\gamma \rightarrow \forall x \psi)$ entails $\forall \tuple v \forall x ((\exists x \gamma) \rightarrow \psi)$ in first order logic, where $\tuple v = \free_T(\gamma) \cup \free_T(\psi)$ (and, therefore, if $x$ is free in $\gamma$ then $x$ is in $\tuple v$).
			\end{enumerate}
		\item[PS-dep] By definition, $=\!\!(\tuple t, t')$ stands for $\indep{\tuple t}{t'}{t'}$; therefore, by rule \textbf{PS-ind} we have that 
			\begin{align*}
				&\forall \tuple v_1 \tuple v_2 ((\gamma(\tuple v_1) \wedge \gamma(\tuple v_2) \wedge \tuple t(\tuple v_1) = \tuple t(\tuple v_2)) \rightarrow\\
				& ~~~\exists \tuple v_3 (\gamma(\tuple v_3) \wedge \tuple t t'(\tuple v_3) = \tuple t t'(\tuple v_1) \wedge \tuple t t'(\tuple v_3) = \tuple t t'(\tuple v_2))) ~|~ \gamma \vdash =\!\!(\tuple t, t').
			\end{align*}
			But the formula
			\[
				\forall \tuple v_1 \tuple v_2 ((\gamma(\tuple v_1) \wedge \gamma(\tuple v_2) \wedge \tuple t(\tuple v_1) = \tuple t(\tuple v_2)) \rightarrow t'(\tuple v_1) = t'(\tuple v_2))
			\]
			entails the premise, and therefore by rule \textbf{PS-ent} we have our conclusion.
	\end{description}
\end{proof}
\begin{theorem}[Soundness]
	Suppose that $\Gamma ~|~ \gamma \vdash \phi$ is provable. Then it is valid. 
\end{theorem}
\begin{proof}
	If $S$ is a provable sequent then there exists a proof $S_1 \ldots S_n S$ for it. Then we go by induction of the length $n$ of this proof:
	\begin{description}
		\item[Base case] Suppose that the proof has length $0$. Then S is an instance of \textbf{PS-lit} or of \textbf{PS-ind}. Suppose first that it is the former, that is, that 
			\[
				S = \forall \tuple v (\gamma \rightarrow \phi) ~|~ \gamma \vdash \phi
			\]
			for some first order $\gamma$ and some first order literal $\phi$, where\\ $\tuple v = \free_T(\gamma) \cup \free_T(\phi)$ and $\phi$ has no parameter variables. Now suppose that $M \models_h \forall \tuple x (\gamma \rightarrow \phi)$; then, by definition, if $s$ is an assignment over team variables such that $M \models_{h\cup s} \gamma$ then $M \models_{s} \phi$. Therefore, by \textbf{ES-lit}, $M \models_{\gamma(s)} \phi$ in entailment semantics, as required.\\

			The case corresponding to \textbf{PS-ind} and \textbf{ES-ind} is entirely similar. 
		\item[Induction case] Let $S_1 S_2 \ldots S_n S$ be our proof. For each $i \leq n$ we have that $S_1 \ldots S_i$ is a valid proof for $S_i$, and hence by induction hypothesis that $S_i$ is valid. Now let us consider which rule $r$ was been used to derive $S$ from $S_1 \ldots S_{n}$: 
			\begin{enumerate}
				\item If $r$ was \textbf{PS-lit} or \textbf{PS-ind} then $(S)$ is a proof for $S$ already, and hence by our base case $S$ is valid; 
				\item If $r$ was \textbf{PS-}$\vee$ then $S$ is $\Gamma_1, \Gamma_2, \forall \tuple v(\gamma \leftrightarrow (\gamma_1 \vee \gamma_2)) ~|~ \gamma \vdash \phi_1 \vee \phi_2$, and there exist two $i, j \leq n$ such that $S_i = (\Gamma_1 ~|~ \gamma_1 \vdash \phi_1)$ and $S_j = (\Gamma_2 ~|~ \gamma_2 \vdash \phi_2)$. By induction hypothesis, these sequents are valid.\\ 
					
					Now suppose that $M \models_h \Gamma_1, \Gamma_2, \forall \tuple v(\gamma \leftrightarrow (\gamma_1 \vee \gamma_2))$. Then, since $M \models_h \Gamma_1$, we have that $M \models_{\gamma_1(h)} \phi_1$, and, analogously, since $M \models_h \Gamma_2$ we have that $M \models_{\gamma_2(h)} \phi_2$. Furthermore, $M \models_h \forall \tuple v(\gamma \leftrightarrow \gamma_1 \vee \gamma_2)$, and therefore by rule \textbf{ES-}$\vee$ we have that $M \models_{\gamma} \phi_1 \vee \phi_2$, as required. 
				\item If $r$ was \textbf{PS-}$\wedge$ then $S_n$ is of the form $\Gamma_1, \Gamma_2 ~|~ \gamma \vdash \phi_1 \wedge \phi_2$ and, by induction hypothesis, $\Gamma_1 ~|~ \gamma \vdash \phi_1$ and $\Gamma_2 ~|~ \gamma \vdash \phi_2$ are valid. Now suppose that $M \models_h \Gamma_1, \Gamma_2$; then $M \models_{\gamma(h)} \phi_1$ and $M \models_{\gamma(h)} \phi_2$, and therefore $M \models_{\gamma(h)} \phi_1 \wedge \phi_2$ by \textbf{ES}-$\wedge$. 
				\item If $r$ was \textbf{PS-}$\exists$ then $S_n$ is of the form $\Gamma, \forall \tuple v (\exists x \gamma' \leftrightarrow \exists x \gamma) ~|~ \gamma \vdash \exists x \phi$, where $\Gamma ~|~ \gamma' \vdash \phi$ is valid by induction hypothesis. Now suppose that $M \models_h \Gamma, \forall \tuple v(\exists x \gamma \leftrightarrow \exists x \gamma')$; then $M \models_{\gamma'(h)} \phi$ and $M \models_h \forall \tuple v (\exists x \gamma \leftrightarrow \exists x \gamma')$, and therefore $M \models_{\gamma(h)} \exists x \phi$ by rule \textbf{ES-}$\exists$.
				\item If $r$ was \textbf{PS-}$\forall$ then $S_n$ is of the form $\Gamma, \forall \tuple v (\gamma' \leftrightarrow \exists x \gamma) ~|~ \gamma \vdash \forall x \phi$, where $\Gamma ~|~ \gamma' \vdash \phi$ is valid by induction hypothesis. Now, suppose that $M \models_h \Gamma, \forall \tuple v (\gamma' \leftrightarrow \exists x \gamma)$. Then $M \models_{\gamma'(h)} \phi$, and furthermore $M \models_h \forall \tuple v(\gamma' \leftrightarrow \exists x \gamma)$. Therefore, by rule \textbf{ES-}$\forall$, $M \models_{\gamma(h)} \forall x \phi$, as required.
				\item If $r$ was \textbf{PS-ent} then $S_n$ is of the form $\Gamma' ~|~ \gamma \vdash \phi$, where $\Gamma ~|~ \gamma \vdash \phi$ is valid by induction hypothesis and where $\bigwedge \Gamma \models \bigwedge \Gamma'$ holds in first order logic. Now suppose that $M \models_h \Gamma'$; then $M \models_h \Gamma$, and hence $M \models_{\gamma(h)} \phi$, as required.

				\item If $r$ was \textbf{PS-depar} then $S_n$ is of the form $\exists p \bigwedge \Gamma ~|~ \gamma \vdash \phi$, where $\Gamma ~|~ \gamma \vdash \phi$ holds by induction hypothesis and where the parameter variable $p$ does not occur free in $\gamma$. Now suppose that $M \models_h \exists p \bigwedge \Gamma$; then there exists an element $m \in \dom(M)$ such that, for $h' = h[m/p]$, $M \models_{h'} \Gamma$. Then $M \models_{\gamma(h')} \phi$; but as $p$ does not occur free in $\gamma$ we then have, by Proposition \ref{propo:non_par_irrel}, that $M \models_{\gamma(h)} \phi$ as required.
				\item If $r$ was \textbf{PS-split} then $S_n$ is of the form $(\bigwedge \Gamma_1) \vee (\bigwedge \Gamma_2) ~|~ \gamma \vdash \phi$, where $\Gamma_1 ~|~ \gamma \vdash \phi$ and $\Gamma_2 ~|~ \gamma \vdash \phi$ by induction hypothesis. Now suppose that $M \models_h (\bigwedge \Gamma_1) \vee (\bigwedge \Gamma_2)$. Then $M \models_h \Gamma_1$ or $M \models_h \Gamma_2$; and in either case, $M \models_{\gamma(h)} \phi$, as required.
			\end{enumerate}
	\end{description}
\end{proof}
In order to prove completeness, we first need a lemma: 
\begin{lemma}
	Suppose that $M \models_{\gamma(h)} \phi$. Then there exists a finite $\Gamma$ such that $\Gamma ~|~ \gamma \vdash \phi$ is provable and such that $M \models_h \Gamma$.
\end{lemma}
\begin{proof}
	The proof is by structural induction on $\phi$.
	\begin{enumerate}
		\item If $\phi$ is a first order literal or an independence atom, this follows immediately from a comparison of \textbf{ES-lit} and \textbf{PS-lit}, and of \textbf{ES-ind} and \textbf{PS-ind}.
		\item If $\phi$ is $\psi_1 \vee \psi_2$ and $M \models_{\gamma(h)} \phi$ then, by definition, there exists an assignment $h'$ extending $h$ and two first order formulas $\gamma_1$, $\gamma_2$ such that $M \models_{\gamma_1(h')} \psi_1$, $M \models_{\gamma_2(h')} \psi_2$ and $M \models_{h'} \forall \tuple v (\gamma \leftrightarrow \gamma_1 \vee \gamma_2)$. Let $\tuple p$ be the tuple of parameters in $\dom(h') \backslash \dom(h)$; now, by induction hypothesis we have that there exist $\Gamma_1$ and $\Gamma_2$ such that $\Gamma_1 ~|~ \gamma_1 \vdash \psi_1$ and $\Gamma_2 ~|~ \gamma_2 \vdash \psi_2$ are provable, and such that furthermore $M \models_{h'} \Gamma_1$ and $M \models_{h'} \Gamma_2$. 

			But then the following is a correct proof:
			\begin{enumerate}
				\item $\Gamma_1 ~|~ \gamma_1 \vdash \psi_1$ (Derived before)
				\item $\Gamma_2 ~|~ \gamma_2 \vdash \psi_2$ (Derived before)
				\item $\Gamma_1, \Gamma_2, \forall \tuple v(\gamma \leftrightarrow \gamma_1 \vee \gamma_2) ~|~ \gamma \vdash \phi$ (\textbf{PS-}$\vee$, (a), (b))
				\item $\exists \tuple p (\bigwedge \Gamma_1 \wedge \bigwedge \Gamma_2 \wedge \forall \tuple v(\gamma \leftrightarrow \gamma_1 \vee \gamma_2)) ~|~ \gamma \vdash \phi$ (\textbf{PS-depar}, (c))\footnote{To be entirely formal, this passage consists of $|\tuple p|$ distinct applications of \textbf{PS-depar}, all of which are correct because none of the parameters in $\tuple p$ appear in $\gamma$.}
			\end{enumerate}
			Finally, $M \models_h \exists \tuple p (\bigwedge \Gamma_1 \wedge \bigwedge \Gamma_2 \wedge \forall \tuple v(\gamma \leftrightarrow \gamma_1 \vee \gamma_2))$, as required, because there exists a tuple of elements $\tuple m$ such that $h[\tuple m/\tuple p] = h'$. 
		\item If $\phi$ is $\psi_1 \wedge \psi_2$ and $M \models_{\gamma(h)} \phi$, then $M \models_{\gamma(h)} \psi_1$ and $M \models_{\gamma(h)} \psi_2$. Then, by induction hypothesis, there exist $\Gamma_1$ and $\Gamma_2$ such that $\Gamma_1 ~|~ \gamma \vdash \psi_1$ and $\Gamma_2 ~|~ \gamma \vdash \psi_2$ are provable and such that $M \models_h \Gamma_1 \Gamma_2$. Then by rule \textbf{PS-}$\wedge$, $\Gamma_1 \Gamma_2 ~|~ \gamma \vdash \psi_1 \wedge \psi_2$, as required. 
		\item If $\phi$ is $\exists x \psi$ and $M \models_{\gamma(h)} \phi$, then there exists a tuple $\tuple p$ of parameter variables not in the domain of $h$, a tuple $\tuple m$ of elements of the model and a formula $\gamma'$ such that, for $h' = h[\tuple m/\tuple p]$, $M \models_{\gamma'(h')} \psi$ and $M \models_{h'} \forall \tuple v (\exists x \gamma' \leftrightarrow \exists x \gamma)$. By induction hypothesis, we then have a $\Gamma'$ such that $\Gamma' ~|~ \gamma' \vdash \psi$ and $M \models_{h'} \Gamma'$. 

			Then the following is a valid proof: 
			\begin{enumerate}
				\item $\Gamma' ~|~ \gamma' \vdash \psi$ (Derived before)
				\item $\Gamma', \forall \tuple v (\exists x \gamma' \leftrightarrow \exists x \gamma) ~|~ \gamma \vdash \exists x \psi$ (\textbf{PS-}$\exists$)
				\item $\exists \tuple p(\bigwedge \Gamma' \wedge \forall \tuple v (\exists x \gamma' \leftrightarrow \exists x \gamma)) ~|~ \gamma \vdash \exists x \psi$ (\textbf{PS-depar})
			\end{enumerate}
			Furthermore, $M \models_h \exists \tuple p (\bigwedge \Gamma' \wedge \forall \tuple v(\exists x \gamma' \leftrightarrow \exists x \gamma))$, as required.
		\item If $\phi$ is $\forall x \psi$ and $M \models_{\gamma(h)} \phi$, then there exists a tuple $\tuple p$ of parameter variables not in the domain of $h$, a tuple $\tuple m$ of elements of the model and a formula $\gamma'$ such that $M \models_{\gamma'(h')} \psi$ and $M \models_{h'} \forall \tuple v(\gamma' \leftrightarrow \exists x \gamma)$, where $h' = h[\tuple m /\tuple p]$. By induction hypothesis, we can then find a $\Gamma'$ such that $\Gamma' ~|~ \gamma' \vdash \psi$ is provable and $M \models_{h'} \Gamma'$. 

		Then the following is a valid proof: 
			\begin{enumerate}
				\item $\Gamma' ~|~ \gamma' \vdash \psi$ (Derived before)
				\item $\Gamma', \forall \tuple v (\gamma' \leftrightarrow \exists x \gamma) ~|~ \gamma \vdash \forall x \psi$ (\textbf{PS-}$\forall$)
				\item $\exists \tuple p(\bigwedge \Gamma' \wedge \forall \tuple v (\gamma' \leftrightarrow \exists x \gamma)) ~|~ \gamma \vdash \forall x \psi$ (\textbf{PS-depar})
			\end{enumerate}
			And, once again, the assignment $h$ satisfies the antecedent of the last sequent, as required. 
	\end{enumerate}
\end{proof}
The completeness of our proof system follows from the above lemma and from the compactness and the L\"owenheim-Skolem theorem for First Order Logic:
\begin{theorem}[Completeness]
	Suppose that $\Gamma ~|~ \gamma \vdash \phi$ is valid, where $\Gamma$ is finite. Then it is provable.
\end{theorem}
\begin{proof}
	Since $\Gamma ~|~ \gamma \vdash \phi$ is valid, for any first order model $M$ over the signature of $\Gamma$, $\gamma$ and $\phi$ and for all $h$ such that $M \models_h \Gamma$ we have that $M \models_{\gamma(h)} \phi$, and hence by the lemma that $M \models_h \Gamma_{M, h}$ for some finite $\Gamma_{M, h}$ such that $\Gamma_{M, h} ~|~ \gamma \vdash \phi$ is provable. \\

	Then consider the first order, countable\footnote{The fact that it is countable follows at once from the fact that it is a first order theory over a countable vocabulary.} theory
	\begin{align*}
	T =  \{\bigwedge \Gamma\} \cup \{\lnot \bigwedge \Gamma_{M, h} :& M \mbox{ is a countable model},\\
	& h \mbox{ is an assignment s.t. } M \models_h \Gamma\}.
	\end{align*}
	This theory is unsatisfiable. Indeed, suppose that $M_0$ is a model that satisfies $\bigwedge \Gamma$ under the assignment $h_0$: then, by the L\"owenheim-Skolem theorem, there exists a countable elementary submodel $(M'_0, h'_0)$ of $(M_0, h_0)$. 

	Now, $M'_0 \models_{h'_0} \Gamma$ and $M'_0$ is countable, and hence by definition $M'_0 \models_{h'_0} \Gamma_{M'_0, h'_0}$. 
	
	But then $M_0 \models_{h_0} \Gamma_{M'_0, h'_0}$ too, and therefore $M_0$ is not a model of $T$.

	By the compactness theorem, this implies that there exists a finite subset $T_0 = \{\lnot \bigwedge \Gamma_{M_1, h_1}, \ldots, \lnot \bigwedge \Gamma_{M_n, h_n}\}$ of $T$ such that $\{\bigwedge \Gamma\} \cup T_0$ is unsatisfiable, that is, such that 
	\[
	\Gamma \models (\bigwedge \Gamma_{M_1, h_1}) \vee \ldots \vee (\bigwedge \Gamma_{M_n, h_n}).
	\]
	
	Now, for each $i$, $\Gamma_{M_i, h_i} ~|~ \gamma \vdash \phi$ can be proved. Therefore, by rule \textbf{PS-split}, we have that $(\bigwedge \Gamma_{M_1, s_1}) \vee \ldots \vee (\bigwedge \Gamma_{M_n, s_n}) ~|~ \gamma \vdash \phi$ is also provable; and finally, by rule \textbf{PS-ent} we can prove that $\Gamma ~|~ \gamma \vdash \phi$, as required.
\end{proof}
Using essentially the same method, it is also possible to prove a ``compactness'' result for our semantics: 
\begin{theorem}
	Suppose that $\Gamma ~|~ \gamma \vdash \phi$ is valid. Then there exists a finite $\Gamma_0 \subseteq \Gamma$ such that $\Gamma_0 ~|~ \gamma \vdash \phi$ is provable (and valid).
\end{theorem}
\begin{proof}
	Let $\kappa = \max(|\Gamma|, \aleph_0)$, and consider the theory 
	\[
		T = \Gamma \cup \{\lnot \bigwedge \Gamma_{M, h} : |M| \leq \kappa, M \models_h \Gamma\}
	\]
	where, as in the previous proof, $\Gamma_{M, h}$ is a finite theory such that $M \models_h \Gamma_{M, h}$ and such that $\Gamma_{M, h} ~|~ \gamma \vdash \phi$ is provable in our system.\\

	Then $T$ is unsatisfiable: indeed, if $T$ had a model then it would have a model $(M, h)$ of cardinality at most $\kappa$, and since that model would satisfy $\Gamma$ it would satisfy $\Gamma_{M, h}$ too, which contradicts our hypothesis. 

	Hence, by the compactness theorem, there exists a finite set\\ $\{\bigwedge \Gamma_{M_1, h_1}, \ldots, \bigwedge \Gamma_{M_n, h_n}\}$ and a finite $\Gamma_0 \subseteq \Gamma$ such that 
	\[
		\Gamma_0 \models \bigwedge \Gamma_{M_1, h_1} \vee \ldots \vee \bigwedge \Gamma_{M_n, h_n}.
	\]
	But by rule \textbf{PS-split}, we have that $\bigwedge \Gamma_{M_1, h_1} \vee \ldots \vee \bigwedge \Gamma_{M_n, h_n} ~|~ \gamma \vdash \phi$ is provable, and hence by rule \textbf{PS-ent} $\Gamma_0 ~|~ \gamma \vdash \phi$ is also provable, as required.
\end{proof}
\section{Adding more teams}
The proof system that we developed in the previous section is, as we saw, sound and complete with respect to its intended semantics. However, this semantics is perhaps quite weak: all that we know is that the teams which correspond to parametrized first order formulas belong in our general models.\\

Rather than adding more and more axioms to our proof system in order to guarantee the existence of more teams, in this section we will attempt to separate our assumptions about team existence from our main proof system. This will allow us to \emph{modulate} our formalism: depending on our needs, we may want to assume the existence of more or of less teams in our general model. 

The natural language for describing assertions about the existence of relations is of course, existential second order logic. The following definitions show how it can be used for our purposes: 
\begin{definition}[Relation Existence Theory]
	A \emph{relation existence theory} $\Theta$ is a set of existential second order sentences of the form $\exists \tuple R \phi(\tuple R)$, where $\phi$ is first order. 
\end{definition}
\begin{definition}[$\Theta$-closed general models]
	Let $(M, \mathcal G)$ be a general model, and let $\Theta$ be a relation existence theory. Then $(M, \mathcal G)$ is $\Theta$-\emph{closed} if and only if for all $\exists \tuple R \phi(\tuple R)$ in $\Theta$ there exists a tuple of teams $\tuple X \in \mathcal G$ such that $M \models \phi[\tuple \rel(\tuple X) / \tuple R]$. 
\end{definition}
\begin{definition}[$\Theta$-valid sequents]
	Let $\Gamma ~|~ \gamma \vdash \phi$ be a sequent and let $\Theta$ be a relation existence theory. Then $\Gamma | \gamma \vdash \phi$ is \emph{valid} if and only if for all $\Theta$-closed models $(M. \mathcal G)$ and all parameter assignments $h$ with domain $\free_P(\Gamma) \cup \free_P(\gamma)$ such that $M \models_h \Gamma$ it holds that 
	\[
		(M, \mathcal G) \models_{\|\gamma\|_h} \phi.
	\]
\end{definition}
Our proof system for $\Theta$-closed general models can then be obtained by adding the following rule to our system: 
\begin{description}
	\item[PS-$\Theta$] If $\Gamma_1(\tuple S), \Gamma_2 ~|~ \gamma \vdash \phi$ is provable, where the relation symbols $\tuple S$ do not occur in $\Gamma_2$, in $\gamma$ or in $\phi$, and $\exists \tuple R \bigwedge \Gamma_1(\tuple R)$ is in $\Theta$ for some $\tuple R$ then $\Gamma_2 ~|~ \gamma \vdash \phi$ is provable.
\end{description}
\begin{theorem}[Soundness]
	Let $\Gamma ~|~ \gamma \vdash \phi$ be a sequent which is provable in our proof system plus \textbf{PS-}$\Theta$. Then it is $\Theta$-valid.
\end{theorem}
\begin{proof}
	The proof is by induction on the length of the proof, and follows very closely the one given already. Hence, we only examine the case in which the last rule used in the proof is \textbf{PS-}$\Theta$. Then, by induction hypothesis, we have that $\Gamma_1(\tuple S), \Gamma ~|~ \gamma \vdash \phi$ is $\Theta$-valid for some $\Gamma_1$ and some $\tuple S$ which does not occur in $\Gamma$, in $\gamma$ or in $\phi$, and moreover $\exists \tuple R \bigwedge \Gamma_1(\tuple R)$ is in $\Theta$. 

	Now, let $(M, \mathcal G)$ be any $\Theta$-closed general model, and let us assume without loss of generality that the relation symbols in $\tuple S$ are not part of its signature. Furthermore, let $h$ be a parameter assignment (with domain $\free(\Gamma) \cup \free(\gamma)$) such that $M \models_h \Gamma$. By definition, there exists a tuple of teams $\tuple X \in \mathcal G$ such that $M \models \bigwedge \Gamma_1[\tuple \rel(\tuple X)/\tuple S]$. Now let $M'$ be $M[\tuple \rel(\tuple X) / \tuple S]$: since $\tuple X$ is in $\mathcal G$, it is not difficult to see that $(M', \mathcal G)$ is a general model. Furthermore, it is $\Theta$-closed, $M' \models \Gamma_1$, and $M' \models_h \Gamma$. Hence, $(M', \mathcal G) \models_{\|\gamma\|_h} \phi$; but since the relation symbols $\tuple S$ do not occur in $\gamma$ or in $\phi$, this implies that $(M, \mathcal G) \models_{\|\gamma\|_h} \phi$. 
\end{proof}
In order to prove completeness, we first need a definition and a simple lemma. 
\begin{definition}[$\Theta^{FO}$]
	Let $\Theta$ be a relation existence theory. Then $\Theta^{FO}$ is the theory $\{\theta_i[\tuple S_i / \tuple R] : \exists \tuple R \theta_i (\tuple R) \in \Theta \}$, where the tuples of symbols $\tuple S_i$ are all disjoint and otherwise unused. 
\end{definition}
\begin{lemma}
	Let $\Theta$ be a relation existence theory and let $M$ be a model such that $M \models \Theta^{FO}$. Then the least general model over it $(M, \mathcal L)$ is $\Theta$-closed. 
\end{lemma}
\begin{proof}
	Consider any $\exists \tuple R \theta(\tuple R) \in \Theta$. Then $M \models \theta(\tuple S_i)$, for some tuple of relation symbols $\tuple S_i$ in the signature of $M$. Then, the teams $\tuple X$ associated to the corresponding relations are in $\mathcal L$, and for these teams we have that $M \models \theta[\tuple \rel(\tuple X) /\tuple R]$, as required. 
\end{proof}
\begin{theorem}[Completeness]
	Suppose that $\Gamma ~|~ \gamma \vdash \phi$ is $\Theta$-valid. Then it is provable in our proof system plus \textbf{PS-}$\Theta$. 
\end{theorem}
\begin{proof}
	Let $M$ be any first order model satisfying $\Theta^{FO}$, where we assume that the relation symbols used in the construction of $\Theta^{FO}$ do not occur in $\Gamma$, in $\gamma$ or in $\phi$. Then, by the lemma, $(M, \mathcal L)$ is $\Theta$-closed, and this implies that, for all assignments $h$ such that $M \models_h \Gamma$, $M \models_{\|\gamma\|_h} \phi$. 

	Therefore, $\Theta^{FO}, \Gamma ~|~ \gamma \vdash \phi$ is valid; and hence, for some finite $\Delta \subseteq \Theta^{FO}$ it holds that $\Delta, \Gamma ~|~ \gamma \vdash \phi$ is provable. Now we can get rid of $\Delta$ through repeated applications of rule \textbf{PS-}$\Theta$ and, therefore, prove that $\Gamma~|~ \gamma \vdash \phi$, as required.
\end{proof}
\paragraph{Acknowledgements}
The author wishes to thank Jouko V\"a\"an\"anen for suggesting this general-models-based approach to the proof theory of Dependence Logic, as well as for a number of useful suggestions and comments. Furthermore, he thankfully acknowledges the support of the EUROCORES LogICCC LINT programme.
\bibliographystyle{plain}
\bibliography{biblio}
\end{document}